\documentclass[a4,11pt]{article}
\usepackage{amsmath,amsfonts,amsthm,latexsym,amssymb,mathrsfs}

\def\H_0{\mathcal{H}_0(T)}

\def\C{{\mathbb C}}

\def\C{{\mathbf C}}

\def\X{{\mathcal X}}
\def\S{\sigma_{SBW_+^-}(T)}
\def\ind{{\textrm{ind}}}
\def\asc{{\textrm{asc}}}
\def\dsc{{\textrm{desc}}}
\def\iso{{\textrm{iso}}}

\def\ind{\textrm{ind}}
\def\X{{\cal X}}
\def\Y{{\cal Y}}
\def\H{{\cal H}}

\def\S{\tau_{AB}}

\def\iso{\textrm{iso}}

\def\asc{ \textrm{asc}}
\def\dsc{ \textrm{dsc}}
\newtheorem{df}{Definition}[section]
\newtheorem{thm}[df]{Theorem}
\newtheorem{pro}[df]{Proposition}

\newtheorem{rema}[df] {Remark}
\newtheorem{lem}[df] {Lemma}
\def\sfstp{{\hskip-1em}{\bf.}{\hskip1em}}

\def\enddemo{\qed \endtrivlist}
\expandafter\let\csname enddemo*\endcsname=\enddemo

\def\qedsymbol{\ifmmode\bgroup\else$\bgroup\aftergroup$\fi
  \vcenter{\hrule\hbox{\vrule
height.5em\kern.5em\vrule}\hrule}\egroup}
\def\qed{\ifmmode\else\unskip\nobreak\fi\quad\qedsymbol}

\title
{ \bf  Generalized Browder's and  Weyl's theorems\\ for left
and right
multiplication operators   \/}

\author { E. $\rm Boasso$, B. P. $\rm Duggal$, I. H. $\rm Jeon$}

\date{   }

\begin{document}

\maketitle \thispagestyle{empty} 

\setlength{\baselineskip}{14pt}

\begin{abstract}\noindent The main objective of this work is to study  generalized
Browder's and Weyl's theorems for the multiplication operators $L_A$ and $R_B$
and for the elementary operator $\S=L_AR_B$.\par
\noindent  \it Keywords: \rm  Banach space,
 left and right multiplication operators, single valued extension property, generalized
 Browder's and  Weyl's theorems.

\indent
\end{abstract}


\section {\sfstp Introduction}\setcounter{df}{0}
\ \indent Let $\X$ and $\Y$ be two Banach spaces and consider two
operators $A\in B(\X)$ and $B\in B(\Y)$. Let $L_A\in B(B(\X))$ and
$R_B\in B(B(\Y))$ be the left and right multiplication operators,
respectively, and denote by $\S=L_AR_B\in B(B(\Y,\X))$ the
elementary operator defined by $A$ and $B$. The main objective of
the present article is to study the generalized Browder's and
Weyl's theorems for $L_A$, $R_B$ and $\S$ (concerning notation
and the main concepts used in this work, see section 2).\par

\indent In section three the (generalized) Browder's, the
(generalized) $a$-Browder's, the Weyl's and the $a$-Weyl's
theorems will be proved to hold for $L_A$ and $R_B$. In addition,
characterizations for $L_A$ and $R_B$ to satisfy generalized
Weyl's and generalized $a$-Weyl's theorem will be presented.
Furthermore, Browder's theorem will be characterized in terms of
the B-Weyl spectrum of multiplication operators. Similar results
will be proved for the $a$-Browder's theorem.\par

\indent In section four the problem of transferring (generalized)
Browder's and (generalized) $a$-Browder's theorem from $A$ and
$B^*$  to $\S$ will be studied. Furthermore, when $A$ and $B^*$
are isoloid operators for which generalized Weyl's (respectively
generalized $a$-Weyl's) theorem holds, necessary and sufficient
conditions for $\S$ to satisfy generalized Weyl's (respectively
generalized $a$-Weyl's) theorem will be given.  \par

\section {\sfstp Notation and terminology}\setcounter{df}{0}
 \hskip.5truecm
\indent This article is concerned with the transmission of Weyl's and Browder's theorem
from bounded and linear maps defined on Banach spaces to the multiplication and
elementary operators induced by them. Weyl's theorem says that the complement in the
spectrum of one kind of essential spectrum is one kind of point spectrum; note that this
statement already breaks into two complementary parts. It is therefore necessary
to briefly recall, on the one hand, the notion of Fredholm operators and some of its generalizations
and, on the other, some subsets of the point spectrum.
In first place some basic notation is introduced.\par

\indent Unless otherwise stated, from now on $\X$ shall denote an
infinite dimensional complex Banach space, $B(\X)$ the algebra of
all bounded linear maps defined on and with values in $\X$ and
$K(\X)$ the closed ideal of compact operators. Given $A\in
B(\X)$, $N(A)$ and $R(A)$ will stand for the null space and the
range of $A$ respectively. Recall that $A\in B(\X)$ is said to be
\it bounded below\rm, if $N(A)=0$ and $R(A)$ is closed. Denote
the \it approximate point spectrum \rm of $A$ by
$\sigma_a(A)=\{\lambda\in \mathbb C \colon A-\lambda \hbox{ is
not bounded below} \}$, where $A-\lambda$ stands for $A-\lambda
I$, $I$ the identity map of $B(\X)$. Let
$\sigma_s(A)=\{\lambda\in \mathbb C \colon R(A-\lambda)\neq \X\}$
denote
 the \it sujectivity spectrum \rm of $A$. Clearly,  $\sigma_a(A)\cup  \sigma_s(A)
=\sigma(A)$, the spectrum of $A$. \par

\indent Recall that $A\in B(\X)$  is said to be  a
\it Fredholm \rm operator if $\alpha(A)=\dim N(A)$ and $\beta(A)=\dim \X/R(A)$
are finite dimensional, in which case its \it index \rm is given by
$$
\ind(A)=\alpha(A)-\beta (A).
$$
If $R(A)$ is closed and $\alpha (A)$ is finite (respectively $\beta (A)$ is finite),
then $A\in  B(\X)$ is said to be  \it upper \rm (respectively \it lower\rm) \it semi-Fredholm\rm,
while if $\alpha (A)$ and $\beta (A)$ are finite and equal, so that the index is zero,
$A$ is said to be \it Weyl \rm operator.\
These classes of operators generate the
Fredholm or essential spectrum and the upper semi-Fredholm, the  lower
semi-Fredholm and the Weyl spectra of $A\in B(\X)$, which will be denoted by
 $\sigma_e(A)$, $\sigma_{SF_+}(A)$,
$\sigma_{SF_-}(A)$ and $\sigma_w(A)$,
 respectively.  It is worth noticing that $\sigma_w(A)=\cap\{ \sigma(A+K)\colon K\in K(\X)\}$  \cite[Corollary 3.41]{A}.
On the other hand, $\Phi(A)$ and $\Phi_+(A)$ will denote the complement
in $\mathbb C$ of the Fredholm spectrum and of the upper semi-Fredholm spectrum of $A$,
respectively.
\par

\indent In addition, the \it Weyl essential approximate point spectrum
\rm of $A\in B(\X)$ is the set $\sigma_{aw}(A)=\cap\{
\sigma_a(A+K)\colon K\in K(\X)\} =\{\lambda\in \sigma_a(A)\colon
A-\lambda\hbox{ is not upper semi-Fredholm or } 0<\ind (A-\lambda) \}$,
see \cite{R2}. \par

\indent In recent years there have been generalizations of the
Fredholm concept. An operator $A\in B(\X)$ will be said to be
\it Berkani Fredholm \rm or \it B-Fredholm\rm, if there exists
$n\in\mathbb N$ for which the range of $R(A^n)$ is closed and the
induced operator $A_n\in B(R(A^n))$ is Fredholm. In a similar way
it is possible to define upper and lower Berkani Fredholm or
B-Fredholm operators. Note that if for some $n\in\mathbb N$,
$A_n\in B(R(A^n))$ is Fredholm, then $A_m\in B(R(A^m))$ is
Fredholm for all $m\ge n$; moreover $\ind (A_n)=\ind (A_m)$, for
all $m\ge n$. Therefore, it makes sense to define the index of
$A$ by $\ind (A)=\ind (A_n)$. Recall that $A$ is said to be \it
Berkani Weyl \rm or \it B-Weyl\rm, if $A$ is  B-Fredholm and
$\ind(A)=0$. Naturally, from these classes of operators, the
B-Fredholm and the B-Weyl spectra of $A\in B(\X)$ can be derived,
which will be denoted by $\sigma_ {BF}(A)$ and $\sigma_ {BW}(A)$,
respectively. In addition, set $\sigma_ {SBF_+^-}(A)= \{\lambda\in
\mathbb C\colon A-\lambda  \hbox{ is not upper semi B-Fredholm or
} 0<\ind (A-\lambda)\}$, see \cite{B6}.\par

\indent In order to state the (generalized) Weyl's theorem, some
subsets of the point spectrum need to be recalled.
First recall that if $\mathcal{K}\subseteq\mathbb C$,
then iso $\mathcal{K}$ is the set of all isolated points of
$\mathcal{K}$ and acc $\mathcal{K}=\mathcal{K}\setminus$ iso $\mathcal{K}$.
\par

\indent  Let $A\in B(\X)$ and denote by $E(A)= \{\lambda\in
iso\hskip.1truecm \sigma(A)\colon 0<\alpha(A-\lambda)\}$ (respectively, by
$E_0(A)= \{\lambda\in E(A)\colon \alpha(A-\lambda)<\infty\}$) the
set of eigenvalues of $A$ which are isolated in the spectrum of
$A$ (respectively, the eigenvalues of finite multiplicity
isolated in $\sigma(A)$). Similarly, define $E^a(A)= \{\lambda\in
iso\hskip.1truecm \sigma_a(A)\colon 0<\alpha(A-\lambda)\}$ (respectively
$E_0^a(A)= \{\lambda\in E^a(A)\colon \alpha(A-\lambda)<\infty\}$)
the set of eigenvalues of $A$ which are isolated in $\sigma_a(A)$
(respectively, the eigenvalues  of finite multiplicity isolated in
$\sigma_a(A)$). \par

\begin{df}Consider a Banach space $\X$ and $A\in B(\X)$. Then it will be said that\par
(i) Weyl's theorem ($Wt$) holds for $A$, if $\sigma_w(A)=\sigma(A)\setminus E_0(A)$,\par
(ii) generalized Weyl's theorem ($gWt$) holds for $A$, if $\sigma_{BW}(A)=\sigma(A)\setminus E(A)$,\par
(iii) $a$-Weyl's theorem ($a$-$Wt$) holds for $A$, if $\sigma_{aw}(A)=\sigma_a(A)\setminus E^a_0(A)$,\par
(iv) generalized $a$-Weyl's theorem ($a$-$gWt$) holds for $A$, if $\sigma_{SBF_+^-}(A)=\sigma_a(A)\setminus E^a(A)$,\par
\end{df}

\indent For information on characterizations and connections amongst the notions
recalled in Definition 2.1, see \cite{B6,R3,D2}.\par

\indent On the other hand, recall that an operator $T\in B(\X)$ is said to
have SVEP, the single--valued extension property,  at a (complex)
point $\lambda_0$, if for every open disc ${\mathcal D}$ centered
at $\lambda_0$ the only analytic function $f:{\mathcal
D}\longrightarrow \X$ satisfying $(T-\lambda)f(\lambda)=0$ is the
function $f\equiv 0$. We say that $T$ has SVEP on a subset
${\mathcal K}$ of the complex plane if it has SVEP at every point
of $\mathcal{K}$. Trivially, every operator $T$ has SVEP at
points of the resolvent $\rho(A)={\C}\setminus \sigma(T)$. Also
$T$ has  SVEP at points $\lambda\in \textrm{iso}\hskip.1truecm \sigma(T)$ and
$\lambda\in \textrm{iso}\hskip.1truecm \sigma_a(T)$.
See \cite[Chapter 2]{A} for
more information on operators with SVEP.\par

 The following
technical lemma will be used in the sequel, often without further
reference. \begin{lem} If $S\in B(\X)$ has SVEP at $\lambda\in
\sigma(S) \setminus \sigma_{{SF}_+}(S)$, then $\lambda \in
iso\hskip.1truecm \sigma_a(S)$ and $asc(S-\lambda)$ is
finite.\end{lem}\begin{demo} See \cite[Theorem
3.23]{A}.\end{demo} It is immediate from the lemma that
$\lambda\in\sigma(S)\setminus\sigma_{ab}(S)$ if and only if
$\lambda\in\sigma_a(S)\setminus\sigma_{{SF}_+}(S)$ and $S$ has
SVEP at $\lambda$ ({\em{cf.}} \cite[Corollary 2.2]{R2}).

\section {\sfstp Operators $L_A$ and $R_A$}\setcounter{df}{0}
\ \indent In this section  Weyl's and Browder's theorems for
 multiplication operators will be studied. First we recall Browder's theorem. To this end, several spectra and
subsets of isolated points need to be considered.\par

\indent The \it ascent \rm  (respectively \it descent\rm ) of
$A\in B(\X)$ is the smallest non-negative integer $a$
(respectively $d$) such that $N(A^a)=N(A^{a+1})$ (respectively
$R(A^d)=R(A^{d+1})$); if such an integer does not exist, then
$asc(A)=\infty$ (respectively $dsc(A)=\infty$). The
operator $A$ will be said to be \it Browder\rm,
if it is Fredholm and its ascent and descent are finite. Then,
the Browder spectrum of $A\in B(\X)$ is the set
$\sigma_b(A)=\{\lambda\in \mathbb C\colon A-\lambda\hbox{ is not Browder}\}$. It is well known that
$$
\sigma_e(A)\subseteq \sigma_w(A)\subseteq \sigma_b(A)= \sigma_e(A)\cup acc\hskip .1truecm\sigma (A).
$$
It is also well known, \cite[Theorem 3.48]{A}, that
$\sigma_b(A)=\cap\{ \sigma(A+K)\colon K\in K(\X), KA=AK\}$.

\indent In addition, the \it Browder essential  approximate point
spectrum \rm of $A\in B(\X)$ is the set
$\sigma_{ab}(A) =\cap\{ \sigma_a(A+K)\colon AK=KA,\hskip.1truecm
K\in K(\X)\} =\{\lambda\in \sigma_a(A)\colon \lambda\in \sigma_{aw}(A)
 \hbox{ or } asc(A-\lambda)=\infty\}$, see
 \cite{R2}. It is clear that
$\sigma_{aw}(A)\subseteq \sigma_{ab}(A)\subseteq \sigma_a(A)$.\par

\indent Moreover, $A$ is said to be a \it B-Browder \rm operator, if
there is some $n\in \mathbb N$ such that $R(A^n)$ is closed, $A_n\in B(R(A^n))$
is Fredholm and $asc(A_n)$ and $dsc(A_n)$ are finite. As
in section 2, the B-Browder spectrum
of $A\in B(\X)$ can be derived; this set will be denoted by
$\sigma_ {BB}(A)$, see
\cite{B6}.\par

\indent On the other hand, recall that a Banach space operator $A\in B(\X)$ is said to be Drazin invertible,
if there exists a necessarily unique $B\in B(\X)$ and some $m\in \mathbb N$ such that
$$
A^m=A^mBA, \hskip.3truecm BAB=B, \hskip.3truecm AB=BA,
$$
see for example \cite{Bo3, D}. If $DR(B(\X))=\{ A\in B(\X)\colon
A\hbox{ is Drazin invertible} \}$, then the Drazin spectrum of
$A\in B(\X)$ is the set $\sigma_{DR}(A)=\{\lambda\in C\colon
A-\lambda\notin DR(B(\X) \}$, see \cite{Bo3}. Note, \cite[Theorem
4.3]{B5}, that $\sigma_{BW}(A)=\cap \{\sigma_{DR}(A+F)\colon F\in
F(\X)\}$, where $F(\X)$ is the ideal of finite range operators
defined on $\X$. Evidently, $\sigma_{BB}(A)=\sigma_{DR}(A)$. What
is more, according to \cite[Theorem 3.2]{BBO} and a duality
argument, $\sigma_{BB}(A)= \sigma_{DR}(A)= \cap
\{\sigma_{DR}(A+F)\colon F\in F(\X), \hskip.1truecm AF=FA\}$. In
particular, $\sigma_{BW}(A)\subseteq \sigma_{DR}(A)\subseteq
\sigma(A)$.
\par

\indent Next denote by $LD(\X)= \{ A\in B(\X)\colon \hbox{  }a=
\asc(A)<\infty\hbox{  and } R(A^{a+1 }) \hbox{ is closed}\}$ the
set of \it left Drazin invertible \rm operators. Then, given
$A\in B(\X)$, the \it left Drazin spectrum \rm of $A$ is the set
$\sigma_{LD}(A)= \{\lambda\in \mathbb C\colon A-\lambda \notin
LD(\X)\}$. Note that according to \cite[Lemma 2.12]{B6}, $\sigma_
{SBF_+^-}(A)\subseteq \sigma_{LD}(A)\subseteq \sigma_a(A)$.\par

\indent Similarly, denote by $RD(\X)= \{ A\in B(\X)\colon \hbox{
}d= dsc(A)<\infty\hbox{  and } R(A^{d}) \hbox{ is closed}\}$ the
set of \it right Drazin invertible \rm operators. Then, given
$A\in B(\X)$, the \it right Drazin spectrum \rm of $A$ is the set
$\sigma_{RD}(A)= \{\lambda\in \mathbb C\colon A-\lambda \notin
RD(\X)\}$. Concerning the left and the right Drazin spectra, see
for example \cite{MM,B6,BBO}.\par

\indent Let $A\in B(\X)$ and denote by $\Pi (A)=
\{\lambda\in \mathbb C\colon
asc(A-\lambda)=dsc(A-\lambda)<\infty\}$ (respectively
$\Pi_0(A)=\{\lambda\in \Pi(A)\colon \alpha(A-\lambda)<\infty\}$)
the set of poles of $A$ (respectively the poles of finite rank of
$A$). Similarly, denote by $\Pi ^a(A)= \{\lambda\in
iso\hskip.1truecm \sigma_a(A)\colon a=asc(A-\lambda)<\infty \hbox{ and
}R(A-\lambda)^{a+1} \hbox{ is closed}\}$ (respectively
$\Pi_0^a(A)= \{\lambda\in\Pi ^a(A)\colon
\alpha(A-\lambda)<\infty\}$) the set of left poles of $A$
(respectively, the left poles of finite rank of $A$).\par

\begin{df}Consider a Banach space $\X$ and $A\in B(\X)$. Then it will be said that\par
(i) Browder's theorem ($Bt$) holds for $A$, if $\sigma_w(A)=\sigma(A)\setminus \Pi_0(A)$,\par
(ii) generalized Browder's theorem ($gBt$) holds for $A$, if $\sigma_{BW}(A)=\sigma(A)\setminus \Pi(A)$,\par
(iii) $a$-Browder's theorem ($a$-$Bt$) holds for $A$, if $\sigma_{aw}(A)=\sigma_a(A)\setminus \Pi_0^a(A)$,\par
(iv) generalized $a$-Browder's theorem ($a$-$gBt$) holds for $A$, if $\sigma_{SBF_+^-}(A)=\sigma_a(A)\setminus \Pi^a(A)$.
\end{df}

\indent Note that necessary and sufficient for $gBt$
(respectively $a$-$gBt$) to hold is that $Bt$ (respectively
$a$-$Bt$) holds, see \cite[Theorem 2.1]{AZ} (respectively
\cite[Theorem 2.3]{AZ}).  In addition,  it is not difficult to
prove that given  $A\in B(\X)$, $\X$ a Banach space, necessary
and sufficient for $A$ to satisfy $Bt$ (respectively $gBt$) is the
fact that $A^*\in B(\X^*)$ satisfies $Bt$ (respectively $gBt$),
where $\X^*$ stands for the dual space of $\X$ and $A^*$ for the adjoint of $A$.\par

\indent Moreover, if $A\in B(\X)$, then
$\sigma(A)\setminus\sigma_b(A)=\Pi_0(T)$ \cite[Lemma 3.4.2]{CPY}.
Consequently, necessary and sufficient for $A$ to satisfy $Bt$
is the identity $\sigma_w(A)=\sigma_b(A)$, equivalently
acc $\sigma(A)\subseteq \sigma_w(A)$. Moreover, \cite[Theorem
2.1]{AZ},
 these conditions are also equivalent to $\sigma_{BW}(A)=\sigma_{DR}(A)$
and hence to acc $\sigma(A)\subseteq \sigma_{BW}(A)$.\par

\indent Furthermore, according to \cite[Corollary 1.3.3]{CPY},
\cite[Corollary 1.3.4]{CPY} and \cite[Corollary 2.2]{R2},
$\sigma_{ab}(A)=\sigma_a(A)\setminus \Pi_0^a(A)$. Therefore, $A$
satisfies $a$-$Bt$ if and only if $\sigma_{aw}(A)=\sigma_{ab}(A)$.
Concerning $a$-$gBt$, since $\sigma_a(A)\setminus
\Pi^a(A)=\sigma_{LD}(A)$, a necessary and sufficient condition for
$A$ to satisfy $a$-$gBt$ is the fact that
$\sigma_{SBF_+^-}(A)=\sigma_{LD}(A)$.\par

\indent For further information on characterizations and connections amongst the notions
recalled in Definition 2.1 and Definition 3.1, see \cite{B6, D2}.\par

\indent In what follows multiplication operators will be studied.\par

\begin{thm} If $A\in B(\X)$ is arbitrary, then Browder's theorem holds for $L_A$ and $R_A$.
Also\par \noindent (i) $\sigma_{BW}(A)\subseteq
\sigma_{BW}(L_A)=\sigma_{DR}(L_A)=\sigma_{DR}(A)$\par \noindent
and\par \noindent (ii) $\sigma_{BW}(A)\subseteq
\sigma_{BW}(R_A)=\sigma_{DR}(R_A)=\sigma_{DR}(A).$\par \noindent
Furthermore, each of the following is equivalent to Browder's
theorem for $A$:\par \noindent (iii)
$\sigma_{BW}(L_A)=\sigma_{BW}(A)$,\par \noindent (iv)
$\sigma_{BW}(R_A)=\sigma_{BW}(A)$.
\end{thm}
\begin{demo}Generally,
$$
\sigma_e(T)\subseteq \sigma_w(T)\subseteq \sigma_b(T)\subseteq \sigma(T),
$$
$T\in B(\X)$, and since  \cite[Corollary 3.4]{Esch} $\sigma_e(L_A)=\sigma_e(R_A)=\sigma (A)=\sigma(R_A)=\sigma(L_A)$
the first assertion is clear. Since $L_A$ satisfies Browder's theorem, it also satisfies the generalized Browder's theorem
(\cite[Theorem 2.1]{AZ}), which is equivalent to the first equality of (i), while the second is   \cite[Theorem 4(iv)]{Bo3},
and now the first inclusion follows. The argument for (ii) is the same.\par
\indent Finally, if   $\sigma_{BW}(L_A)\subseteq \sigma_{BW}(A)$, then $\sigma_b (A)\subseteq \sigma_w (A)$.\par
\indent In fact, if $\sigma_{BW}(L_A)\subseteq \sigma_{BW}(A)$, then according to statement (i),
 $\sigma_{BW}(A)=  \sigma_{DR}(A)$. Therefore, according to \cite[Theorem 2.1]{AZ}, Browder's
theorem holds for $A$. On the other hand, if $A$ satisfies Browder's theorem, then
according to  \cite[Theorem 2.1]{AZ},  $\sigma_{BW}(A)=  \sigma_{DR}(A)$.
Consequently, according to (i), statement (iii) holds.\par

\indent A similar argument proves that Browder's theorem holds for $A$ if and only if statement (iv) holds.
\end{demo}

\indent Trivially, Weyl's theorem holds for all operators $A\in
B(\X)$ when $\X$ is finite dimensional. The following theorem
says that this remains true for  operators $L_A$ and $R_A$ for
every $A\in B(\X)$.
\par

\begin{thm} If  $A\in B(\X)$, then  \par
\noindent (i) $E_0(L_A)=E_0(R_A)=\emptyset$,\par
\noindent and hence Weyl's theorem holds for both $L_A$ and $R_A$.
\end{thm}
\begin{demo} It is easy to see (\cite[Theorem 4]{HK}) that each of $N(L_A)$ and $N(R_A)$ are either zero or infinite dimensional,
giving (i). This together with the first part of Theorem 3.2 finishes the proof.
\end{demo}

\indent In order to study the generalized Weyl's theorem ($gWt$),
some preparation is needed.\par

\begin{lem} If $A\in B(\X)$, then $E(L_A)=E(A)$
and $E(R_A)=E(A^*)$.
\end{lem}
\begin{demo} In first place, it is clear the the isolated points of $\sigma(A)$ and
of $\sigma(L_A)$ coincide. Let $\lambda$ be an eigenvalue of $A$ and consider
$v\in \X$ such that $(A-\lambda )(v)=0\neq v$. Let $H$ be a closed
subspace of $\X$ such that $H\oplus <v>=\X$, and construct $P\in B(\X)$
such that $P\mid_H=0$ and $P(v)=v$. Then, $P\neq 0$ and $(L_A-\lambda)(P)=0$.\par
\indent On the other hand, if there exists $S\in B(\X)$, $S\neq 0$, such that $(L_A-\lambda )(S)=0$,
then there is $y\in \X$, $y \neq 0$, such that $S(y)\neq 0$ and
$(A-\lambda )(S(y))=0$.\par

\indent Concerning the last statement, using adjoint operators it
is clear that $E(R_A)\subseteq E(L_{A^*})=E(A^*)$. On the other
hand, if $\lambda\in E(A^*)$, then there is $f\in {\X}^*$,
$f\neq0$, such that $f(A-\lambda)=0$. Let $v\in \X$, $v\neq 0$,
and define $T\in B(\X)$ as follows: $T(y)=f(y)v$, $y\in \X$.
Then, $T$ is a bounded and linear map, $T\neq 0$, and $T\in N(R_
{A-\lambda})$. Since $iso\hskip.1truecm \sigma (A^*)=iso\hskip.1truecm\sigma( R_A)$,
$E(R_A)=E(A^*)$.
\end{demo}

\begin{thm} If $A\in B(\X)$, then
the following statements are equivalent.\par \noindent
\hskip.2truecm (i) $A$ satisfies $gWt$.\par \noindent (ii) $L_A$
satisfies $gWt$ and  $A$ satisfies $Bt$.\par
\end{thm}
\begin{demo} If the second statement holds, then according to Theorem 3.2(iii), Lemma 3.4
and the identity $\sigma (A)=\sigma (L_A)$,
$$
\sigma_{BW}(A)=\sigma_{BW}(L_A)=\sigma (A)\setminus E(A).
$$
\indent On the other hand, if the first statement holds, then
(since $gWt\Longrightarrow Wt\Longrightarrow Bt$, see \cite{B6})  Browder's
theorem holds for $A$. Therefore (according again to Theorem
3.2(iii), Lemma 3.4 and the fact that $\sigma (A)=\sigma (L_A)$),
$\sigma_{BW}(L_A)=\sigma_{BW}(A)=\sigma (L_A)\setminus
E(L_A).$\end{demo}

\begin{thm} If $A\in B(\X)$, then
the following statements are equivalent.\par \noindent
\hskip.2truecm (i)   $A^*$ satisfies $gWt$.\par \noindent (ii)
$R_A$ satisfies $gWt$ and  $A$ satisfies $Bt$.\par
\end{thm}
\begin{demo}If $A^*$ satisfies $gWt$, then both $A$ and $A^*$ satisfy
$Bt$ (in particular, see Theorem 3.2(iv)),
$\sigma_{BW}(R_A)=\sigma_{BW}(A)$)
 and $\sigma_{BW}(A^*)= \sigma(A^*)\setminus E(A^*)$. Thus
 $\sigma_{BW}(R_A)=\sigma_{BW}(A)=\sigma_{BW}(A^*)$. Applying Lemma 3.4  and the identity $\sigma(A^*)
=\sigma (R_A)$, it follows that
$$
\sigma_{BW}(R_A)= \sigma(R_A)\setminus E(R_A),
$$
equivalently, $R_A$ satisfies $gWt$.\par

\indent On the other hand, if the second statement holds, then
 Theorem 3.2(iv) and  the
fact that both $A$ and $A^*$ satisfy ($Bt$, hence equivalently)
$gBt$ imply that
$\sigma_{BW}(R_A)=\sigma_{BW}(A)=\sigma_{BW}(A^*)$. Therefore, by
Lemma 3.4 and the identity $\sigma (A^*)=\sigma (R_A)$, $A^*$
satisfies $gWt$.
\end{demo}

\indent Next $a$-Weyl's and $a$-Browder's theorems will be
studied. For a subset $K$ of the set of complex numbers, let
$K^*$ denote its complex conjugate. \par

\begin{thm} If $A\in B(\X)$, then
the following statements hold.\par \noindent (i) The operators
$L_A$ and $R_A\in B(B(\X))$ satisfy $a$-$Bt$.\par \noindent (ii)
$A$ satisfies $a$-$Bt$ and $\Pi^a(A)=\Pi^a(L_A)$ if and only if
$\sigma_{SBF_+^-}(A)=\sigma_{SBF_+^-}(L_A)$.\par \noindent (iii)
If $\X$ is an infinite dimensional complex Hilbert space, then a
necessary and sufficient for $A^*$ to satisfy $a$-$Bt$ is the
fact that $
\sigma_{SBF_+^-}(A^*)=(\sigma_{SBF_+^-}(R_A))^*$.
\end{thm}
\begin{demo} $(i)$ Recall that given an operator $A\in B(\X)$, $\sigma_{SF_+}(A)
\subseteq\sigma_{aw}(A)\subseteq \sigma_{ab}(A)\subseteq
\sigma_a(A)$, and
$\sigma_{SF_+}(L_A)\subseteq\sigma_a(L_A)=\sigma_a(A)$,
\cite[Lemma]{DR}. Thus, since   $\sigma_a(A)\subseteq
\sigma_{SF_+}(L_A)$ \cite [Proposition 6.2]{Bo1}, it follows that
 $$\sigma_{SF_+}(L_A)=\sigma_{aw}(L_A)=\sigma_{ab}(L_A)=\sigma_a(L_A)=\sigma_a(A).$$
 A similar argument proves that  $$\sigma_{SF_+}(R_A)=\sigma_{aw}(R_A)=\sigma_{ab}(R_A)=
\sigma_a(R_A)=\sigma_s(A).$$ In particular,  $\sigma_{aw}(L_A)=
\sigma_{ab}(L_A)$ and $\sigma_{aw}(R_A)= \sigma_{ab}(R_A)$, i.e.,
 $L_A$ and $R_A$ satisfy $a$-$Bt$.\par

$(ii)$  If $\sigma_{SBF_+^-}(A)=\sigma_{SBF_+^-}(L_A)$, then
$\sigma_a(A)\setminus\sigma_{SBF_+^-}(A)=\sigma_a(L_A)\setminus\sigma_{SBF_+^-}(L_A)$.
Since $L_A$ satisfies $a$-$Bt$, see (i), $L_A$ satisfies
$a$-$gBt$, i.e.,
$\sigma_a(L_A)\setminus\sigma_{SBF_+^-}(L_A)=\Pi^a(L_A)$. Thus,
since $\Pi^a(L_A)\subseteq\Pi^a(A)$ \cite[Theorem 8(ii)]{Bo3}
and $\Pi^a(A)\subseteq\sigma_a(A)\setminus\sigma_{SBF_+^-}(A)$,
$\Pi^a(A)\subseteq\sigma_a(A)\setminus\sigma_{SBF_+^-}(A)=\Pi^a(L_A)\subseteq\Pi^a(A)$,
i.e., $A$ satisfies $a$-$gBt$ (hence also, $a$-$Bt$) and
$\Pi^a(A)=\Pi^a(L_A)$. Conversely, if $A$ satisfies $a$-$Bt$
(hence also, $a$-$gBt$) and $\Pi^a(A)=\Pi^a(L_A)$, then
$\sigma_a(A)\setminus\sigma_{SBF_+^-}(A)=\Pi^a(A)=\Pi^a(L_A)$.
Furthermore, since $L_A$ satisfies $a$-$gBt$ (by $(i)$ above),
$\sigma_a(L_A)\setminus\sigma_{SBF_+^-}(L_A)=\sigma_a(A)\setminus
\sigma_{SBF_+^-}(L_A)=\Pi^a(L_A)=\sigma_a(A)\setminus\sigma_{SBF_+^-}(A)$,
which implies that $\sigma_{SBF_+^-}(A)=\sigma_{SBF_+^-}(L_A)$.
\par

\indent $(iii)$ Recall, \cite[Proposition 3.10]{D1}, that a
Banach space operator $T$ satisfies $a$-$Bt$ (equivalently,
$a$-$gBt$) if and only if $\sigma_{SBF^-_+}(T)=\sigma_{LD}(T)$.
Since $R_A$ satisfies $a$-$Bt$ (by $(i)$ above),
$\sigma_{SBF_+^-}(R_A)=\sigma_{LD}(R_A)$. Now recall from \cite
[Theorem 9]{Bo3}  and  \cite[ page 139]{MM} that if $\X$ is a complex Hilbert space, then
$\sigma_{LD}(A^*)=(\sigma_{RD}(A))^*=(\sigma_{LD}(R_A))^*$.
Hence $A^*$ satisfies $a$-$Bt$ if and only if
$\sigma_{SBF^-_+}(A^*)=(\sigma_{SBF_+^-}(R_A))^*$.
\end{demo}

\begin{thm} If $A\in B(\X)$, then
$a$-$Wt$ holds for both $L_A$ and $R_A$.
\end{thm}
\begin{demo} According to the proof of Theorem 3.7, it is enough to prove that $E_0^a(L_A)=\emptyset=E_0^a(R_A)$.
However, these identities can be deduced from \cite[Theorem 4]{HK}.
\end{demo}

\begin{thm} If $A\in B(\X)$, then the following statements
are equivalent.\par \noindent (i)  $A$ satisfies $a$-$gWt$ and
$\Pi^a(A)=\Pi^a(L_A)$.\par \noindent (ii)The operator $L_A$
satisfies $a$-$gWt$ and $A$ satisfies $a$-$Bt$.
\end{thm}
\begin{demo}Recall from  \cite[Theorem 8(ii)]{Bo3} that
$\Pi^a(L_A)\subseteq \Pi^a (A)$; since  $E^a(A)=E^a(L_A)$ by
Lemma 3.4, and since $\Pi^a(T)\subseteq E^a(T)$ for every
operator $T$, it follows that $\Pi^a(L_A)\subseteq \Pi^a
(A)\subseteq E^a(A)=E^a(L_A)$. Suppose now that $(ii)$ is
satisfied. Then $L_A$ satisfies generalized $a$-Weyl's theorem
implies that $\Pi^a(L_A)=E^a(L_A)$. Consequently,
$\Pi^a(L_A)=\Pi^a(A)= E^a(A)$, and it follows from the hypothesis
 ``$A$ satisfies $a$-$Bt$" (equivalently, $A$ satisfies
  $a$-$gBt$ $\Longleftrightarrow \sigma_a(A)\setminus\sigma_ {SBF_+^-}(A)=
 \Pi^a(A)$) that $\sigma_a(A)\setminus\sigma_ {SBF_+^-}(A)=E^a(A)$,
 i.e., $A$ satisfies $a$-$gWt$.\par \indent On the
other hand, if the first statement holds, then  $A$ satisfies
($a$-$Bt$, equivalently)  $a$-$gBt$, equivalently
$\sigma_a(A)\setminus\sigma_{SBF_+^-}(A)=\Pi^a(A)$, and
$\Pi^a(A)=E^a(A)$. In view of the  hypothesis
$\Pi^a(A)=\Pi^a(L_A)$, it now follows from Theorem 3.7(ii)  that
$\sigma_{SBF_+^-}(A)= \sigma_{SBF_+^-}(L_A)$. Consequently, since
$\sigma_a(L_A)=\sigma_a(A)$ and $L_A$ satisfies $a$-$Bt$ (by
Theorem 3.7(i)),
$\sigma_a(L_a)\setminus\sigma_{SBF_+^-}(L_A)=E^a(L_a)$, i.e.,
$L_A$ satisfies  $a$-$gWt$.
\end{demo}

\begin{rema}\rm Note that if $\X$ is a Banach space and $A\in B(\X)$,
then $L_A$ satisfies $a$-$gBt$ and the following implications
hold: $$\lambda\notin\sigma_{SBF_+^-}(L_A)\Longleftrightarrow
\lambda\in\Pi^a(L_A)\subseteq\Pi^a(A)\Longrightarrow
\lambda\notin\sigma_{SBF_+^-}(A).$$ (Here the final implication
follows from \cite[Lemma 3.1]{D1}.) Thus
$\sigma_{SBF_+^-}(A)\subseteq\sigma_{SBF_+^-}(L_A)$: the
condition in Theorem 3.7(ii) is equivalent to
$\sigma_{SBF_+^-}(L_A)\subseteq \sigma_{SBF_+^-}(A)$.\par \indent
Note also that if $X$ is a Hilbert space, then (according to
\cite[Lemma]{DR} and \cite[Theorem 9(iii)]{Bo3}) $\Pi^a (A)=\Pi^a
(L_A)$ for all $A\in B(\X)$.\par
\end{rema}

\begin{thm}Let $\X$ be an infinite dimensional complex Hilbert space and
consider an operator $A\in B(\X)$. Then, the following statements
are equivalent.\par \noindent (i)  $a$-$gWt$ holds for $A^*$.\par
\noindent (ii)The operator $R_A$ satisfies $a$-$gWt$ and $A^*$
satisfies $a$-$Bt$.
\end{thm}
\begin{demo} If $\X$ is a complex Hilbert space, then
$E^a(A^*)=(E^a(R_A))^*$ (see Lemma 3.4 and \cite[Lemma]{DR}).
Recall also that $\sigma_a(A^*)=(\sigma_a(R_A))^*$, $A^*$
satisfies $a$-$gBt$ if and only if
$\sigma_{SBF_+^-}(A^*)=(\sigma_{SBF_+^-}(R_A))^*$ (Theorem
3.7(iii)), and $\Pi^a(A^*)=(\Pi^a(R_A))^*$ \cite[Theorem
9]{Bo3}. Hence
\begin{eqnarray*}A^* \hspace{1mm}\mbox{satisfies}\hspace{1mm}
a-gWt & \Longleftrightarrow &
\sigma_a(A^*)\setminus\sigma_{SBF_+^-}(A^*)=\Pi^a(A^*)=E^a(A^*)\\&
\Longleftrightarrow & A^* \hspace{1mm}\mbox{satisfies}\hspace{1mm}
a-gBt\hspace{1mm}\mbox{and}\hspace{2mm}\Pi^a(A^*)=E^a(A^*)\\&
\Longleftrightarrow & A^* \hspace{1mm}\mbox{satisfies}\hspace{1mm}
a-Bt\hspace{1mm}\mbox{and}\hspace{1mm}\sigma_a(A^*)\setminus\sigma_{SBF_+^-}(A^*)=E^a(A^*)\\
& \Longleftrightarrow & A^*
\hspace{1mm}\mbox{satisfies}\hspace{1mm}
a-gBt\hspace{1mm}\mbox{and}\hspace{1mm}(\sigma_a(R_A))^*\setminus (\sigma_{SBF_+^-}(R_A))^*\\
& \hspace{2cm}&=  (E^a(R_A))^*
(=(\Pi^a(R_A))^*=\Pi^a(A^*)=E^a(A^*)) \\
& \Longleftrightarrow & A^*
\hspace{1mm}\mbox{satisfies}\hspace{1mm}
a-gBt\hspace{1mm}\mbox{and}\hspace{1mm}
\sigma_a(R_A)\setminus\sigma_{SBF_+^-}(R_A)\\
& \hspace{2cm}&=E^a(R_A).\\\end{eqnarray*}
Thus $(i)\Longleftrightarrow (ii).$\end{demo}

\section {\sfstp The operator $\S=L_AR_B$}\setcounter{df}{0}
\ \indent In the following $A$ shall denote an operator in
$B(\X)$, $B$ an operator in $B(\Y)$ and $\S=L_AR_B\in B(B(\Y,\X))$
the operator $\S(X)=L_AR_B(X)=AXB$.  Recall from \cite[Corollary
3.4]{Esch} that
$$\sigma(\S)=\sigma(A)\sigma(B)
\hspace{2mm}\mbox{and}\hspace{2mm}\sigma_e(S)=\sigma(A)\sigma_e(B)\cup\sigma_e(A)\sigma(B).$$
Evidently, $\iso\hskip.1truecm \sigma(\S)\subseteq\iso\hskip.1truecm \sigma(A)\iso\hskip.1truecm \sigma(B)$.
A bit more work, \cite[Theorem 4]{HK}, shows that the point
spectrum $\sigma_p(\S)$ of $\S$ satisfies the inclusion $\sigma_p(A)\sigma_p(B^*)\subseteq\sigma_p(\S)$.
However, in order to compute some other spectra that will be relevant for the
present article, some
preliminary definitions should be recalled.\par

\begin{rema}\rm Let $\X$ and $\Y$ be two Banach spaces and $A\in B(\X)$ and $B\in B(\Y)$.
Consider $T$ the two-tuple of commuting operators $T=(L_A,R_B)$,
$L_A$ and $R_B\in B(B(\Y,\X))$. Recall that the \it approximate
joint point spectrum \rm of $T$ and the \it upper semi-Fredholm
joint spectrum \rm of $T$ are the sets
 \begin{align*}
&\sigma_{\pi}(T)=\{ (\mu,\nu)\in \mathbb C^2\colon V(A-\mu,
B-\nu) \hbox{ is not bounded below}\}
 \hbox{ and}&\\
&\sigma_{\Phi_+}(T)=\{ (\mu,\nu)\in \mathbb C^2\colon V(A-\mu, B-\nu) \hbox{ is not upper semi-Fredholm}\},\\
\end{align*}
respectively, where $V(A-\mu, B-\nu)\colon B(\Y,\X)\to
B(\Y,\X)\times B(\Y,\X)$, $V(A-\mu, B-\nu)(S)= (L_{A-\mu}S,
R_{B-\nu}S)= ((A-\mu)S, S(B-\nu))$. Concerning the properties of
these joint spectra, see for example \cite{BHW,Sl,Bo1}.
\end{rema}
The following technical lemma, which will be useful in what
follows, says that in considering points
$\lambda\in\sigma_a(\S)\setminus\sigma_{SF_+}(\S)$, or
$\lambda\in\sigma(\S)\setminus\sigma_e(\S)$, it will suffice to
consider points $\lambda\neq0$.
\begin{lem}\label{lem11}If $A$, $B$ and $\S$ are as above,
then $0\notin\sigma_a(\S)\setminus\sigma_{{SF}_+}(\S)$.\end{lem}\begin{demo}
Suppose to the contrary that
$0\in\sigma_a(\S)\setminus\sigma_{{SF}_+}(\S)$. Then $\S$ is upper
semi--Fredholm, so that $0< \alpha(\S)<\infty$. If $A$ and $B^*$
are injective, then $\S$ is injective. Hence either $0<\alpha(A)$
or $0<\alpha(B^*)$. As in the proof of \cite[Theorem 4]{HK}, let $f\odot x\in B(\Y,\X)$ denote the rank one
operator $(f\odot x)y=f(y)x$ induced by $x\in\X$ and $f\in{\Y}^*$.
Note that $\S(f\odot x)=A(x)\odot B^*(f)$. Then
$$\{\X\odot N({B^*})\}\cup\{N(A)\odot {\Y}^*\}\subseteq
N(\S),$$ which implies that $\alpha(\S)=\infty$. This being a
contradiction, we must conclude that
$0\notin\sigma_a(\S)\setminus\sigma_{{SF}_+}(\S)$.\end{demo}
\indent It is apparent from the argument above that
$0\notin\sigma(\S)\setminus\sigma_e(\S)$.
\begin{pro}If $A$, $B$ and $\S$ are as above, then the following satements hold.\par
 \begin{align*}
&(i)  & &\sigma_a(\S)=\sigma_a(A)\sigma_a(B^*), \hbox{ in particular }iso\hskip.1truecm  \sigma_a(\S)\subseteq iso\hskip.1truecm \sigma_a(A)iso\hskip.1truecm \sigma_a(B^*),&\\
& (ii)& &\sigma_{{SF}_+}(\S)=\sigma_a(A)\sigma_{{SF}_+}(B^*)\cup\sigma_{{SF}_+}(A)\sigma_a(B^*),&\\
& (iii) & &\sigma_b(\S)=\sigma(A)\sigma_b(A)\cup\sigma_b(A)\sigma(B),&\\
& (iv)  & &\sigma_{ab}(\S)=\sigma_a(A)\sigma_{ab}(B^*)\cup\sigma_{ab}(A)\sigma_a(B^*).\\
\end{align*}
\end{pro}
\begin{demo} Concerning statement $(i)$, recall that according to \cite[Proposition 6.1(ii)]{Bo1},
$$
\sigma_a(A)\times \sigma_a(B^*)\subseteq \sigma_{\pi}(L_A,R_B).
$$
However, according to \cite[Lemma]{DR} and the spectral mapping theorem applied to the
maps $P_1$, $P_2\colon \mathbb C^2\to\mathbb C$, the projections on the first and
the second coordinate respectively, \cite[Theorem 2.9]{Sl},
 $$
\sigma_{\pi}(L_A,R_B)=\sigma_a(A)\times \sigma_a(B^*),
$$
which, applying  the spectral mapping theorem to the polynomial mapping
$P\colon \mathbb C^2\to\mathbb C$, $P(x,y)=xy$, \cite[Theorem 2.9]{Sl}, implies that
$\sigma_a(\S)=\sigma_a(A)\sigma_a(B^*)$. The remaining inclusion of the
first statement is clear.\par

\indent To prove  statement $(ii)$, recall that according to
\cite[Proposition 6.2(ii)]{Bo1},
$$
\sigma_{SF_+}(A)\times \sigma_a(B^*)\cup
\sigma_a(A)\times\sigma_{SF_+}(B^*)\subseteq
\sigma_{\Phi_{+}}(L_A,R_B).
$$
However, a direct calculation using the definition of  $V(A-\mu,
B-\nu)\colon B(\Y,\X)\to B(\Y,\X)\times B(\Y,\X)$ proves that the
last inclusion is an equality, which, applying  the spectral
mapping theorem to the polynomial mapping $P\colon \mathbb
C^2\to\mathbb C$, $P(x,y)=xy$, \cite[Theorem 7]{BHW}, implies
statement $(ii)$.\par

\indent Next we consider statement $(iii)$. Let
$\lambda\in\sigma(\S)\setminus\sigma_b(\S)$. Since
$\sigma_b(\S)=\sigma_e(\S)\cup acc\hskip.1truecm \sigma(\S)$, $\lambda\in
iso\hskip.1truecm \sigma(\S)\subseteq iso\hskip.1truecm \sigma(A)iso\hskip.1truecm \sigma(B)$ and, according
to \cite[Corollary 3.4]{Esch}, $\lambda\notin
\sigma_e(\S)=\sigma(A)\sigma_e(B)\cup\sigma_e(A)\sigma(B)$.
Therefore, there are $\mu$ and $\nu$ such that $\lambda=\mu \nu$,
$\mu\in iso\hskip.1truecm  \sigma(A)$, $\nu\in iso \hskip.1truecm \sigma(B)$,
$\mu\in\sigma(A)\setminus\sigma_e (A)$ and
$\nu\in\sigma(B)\setminus\sigma_e(B)$. In particular,
$\mu\in\sigma(A)\setminus\sigma_b (A)$ and
$\nu\in\sigma(B)\setminus\sigma_b(B)$. Therefore, $\lambda\notin
(\sigma(A)\sigma_b(B)\cup\sigma_b(A)\sigma(B))$. Hence,
$\sigma(A)\sigma_b(B)\cup\sigma_b(A)\sigma(B)\subseteq
\sigma_b(\S)$.\par

\indent On the other hand, if $\lambda \in
\sigma(\S)\setminus(\sigma(A)\sigma_b(B)\cup\sigma_b(A)\sigma(B))$,
then, according to \cite[Corollary 3.4]{Esch},  $\lambda\notin
\sigma_e(\S)$, and for every $\mu\in\sigma (A)$ and $\nu\in\sigma
(B)$ such that $\lambda=\mu\nu$, $\mu\in
\sigma(A)\setminus\sigma_b(A)$ and
$\nu\in\sigma(B)\setminus\sigma_b(B)$. However, since for the
factorization of $\lambda=\mu\nu$ with $\mu$ and $\nu$ as before,
$\mu\in iso\hskip.1truecm  \sigma (A)$ and $\nu\in iso\hskip.1truecm   \sigma (B)$,
 $\lambda=\mu\nu\in iso \hskip.1truecm  \sigma(\S)$ \cite[Theorem 6]{HK}.
Consequently, $\lambda\in\sigma(\S)\setminus\sigma_b(\S)$. Hence,
$\sigma_b(\S)\subseteq
\sigma(A)\sigma_b(B)\cup\sigma_b(A)\sigma(B)$.\par

\indent To prove  statement $(iv)$, start by recalling  that
 given an operator $T$
defined on the Banach space $\X$, a necessary and sufficient for
$\lambda\in \sigma_a(T)\setminus\sigma_{ab}(T)$ is that
$\lambda\in iso\hskip.1truecm  \sigma_a(T)$, $\lambda\notin\sigma_{SF_+}(T)$ and
$asc(T-\lambda)$ is finite  \cite[Corollary 2.2]{R2}.  \par

\indent Suppose that $\lambda\in
\sigma_a(\S)\setminus\sigma_{ab}(\S)$. Then  statement $(i)$
(together with \cite[Corollary 2.2]{R2}) implies the existence of
$\mu\in \sigma_a(A)$ and $\nu\in\sigma_a(B^*)$ such that
$\lambda=\mu\nu$,  $\mu\in iso\hskip.1truecm  \sigma_a(A)$ and $\nu\in
iso\hskip.1truecm  \sigma_a(B^*)$. Since $\sigma_{SF_+}(\S)\subseteq
\sigma_{ab}(\S)$, statement $(ii)$ (together with \cite[Corollary
2.2]{R2}) implies that
$\mu\in\sigma_a(A)\setminus\sigma_{SF_+}(A)$ and
$\nu\in\sigma_a(B^*)\setminus\sigma_{SF_+}(B^*)$. But then
$asc(A-\mu)$ and $asc(B^*-\nu)$ are finite (see Lemma 2.2). Hence
 $\mu\in \sigma_a(A)\setminus\sigma_{ab}(A)$ and $\nu\in
\sigma_a(B^*)\setminus\sigma_{ab}(B^*)$. Thus, $\lambda\in
\sigma_a(\S)\setminus
(\sigma_a(A)\sigma_{ab}(B^*)\cup\sigma_{ab}(A)\sigma_a(B^*))$;
hence
$\sigma_a(A)\sigma_{ab}(B^*)\cup\sigma_{ab}(A)\sigma_a(B^*)\subseteq
\sigma_{ab}(\S)$.\par

\indent On the other hand, consider $\lambda\in
\sigma_a(\S)\setminus
(\sigma_a(A)\sigma_{ab}(B^*)\cup\sigma_{ab}(A)\sigma_a(B^*))$.
Then, according to  statement $(ii)$, $\lambda\notin
\sigma_{SF_+}(\S)$ and for every $\mu\in\sigma_a (A)$ and
$\nu\in\sigma_a (B^*)$ such that $\lambda=\mu\nu$, $\mu\in
\sigma_a(A)\setminus\sigma_{ab}(A)$ and $\nu\in\sigma_a(B^*)
\setminus\sigma_{ab}(B^*)$. In particular, $\mu\in
iso\hskip.1truecm  \sigma_a(A)$ and $\nu\in iso \hskip.1truecm  \sigma_a(B^*)$. This, according
to  statement $(i)$ and \cite[Theorem 6]{HK}, implies that
$\lambda\in iso\hskip.1truecm  \sigma_a(\S)$. Additionally, see Lemma 2.2,
$asc(\S-\lambda)$ is finite. Consequently, $\lambda\in
\sigma_a(\S)\setminus \sigma_{ab}(\S)$; hence
$\sigma_{ab}(\S)\subseteq
\sigma_a(A)\sigma_{ab}(B^*)\cup\sigma_{ab}(A)\sigma_a(B^*)$.
\end{demo}

\indent Before studying the generalized Browder's theorem for the
operator $\S$, some results must be considered.\par

 Let
$(0\neq)\lambda\in\sigma(\S)\setminus\sigma_e(\S)$. Let \vskip6pt
$$
E=\{(\mu_i,\nu_i)_{i=1}^p\in\sigma
(A)\sigma(B)\!:\,\mu_i\nu_i=\lambda\}.
$$
\vskip4pt \noindent Then $E$ is a finite set, see \cite[Lemma 4.1]{Esch}. Furthermore,
\cite[Theorem 4.2]{Esch}, \vskip4pt
\begin{description}
\item{\rm$\;$(i)$\:$}
if $n>1$, then $\mu_i\in\iso \hskip.1truecm  \sigma(A)$, for $1\le i\le n$,
\vskip4pt \item{\rm$\,$(ii)$\,$} if $p>n$, then
$\nu_i\in\iso\hskip.1truecm  \sigma(B)$, for $n+1\le i\le p$, \vskip4pt
\item{\rm(iii)}
$\ind(\S-\lambda)=\sum_{j=n+1}^p\ind(A-\mu_j)\dim H_0(B-\nu_j)$
\vskip4pt\hskip70pt $-\sum_{j=1}^n\ind (B-\nu_j)\dim
H_0(A-\mu_j)$.
\end{description}
 Here $H_0(A-\mu_i)$, similarly $H_0(B-\nu_i)$,
denotes the quasinilpotent part
$H_0(A-\mu_i)=\{x\in\X:\lim_{n\longrightarrow\infty}{||(A-\mu_i)^nx||^{\frac{1}{n}}=0}\}$
of the operator $A-\mu_i$, $\mu_i\in\iso\hskip.1truecm  \sigma(A)$. Clearly,
$\dim{H_0(A-\mu_i)}$ ($1\leq i\leq n$), similarly
$\dim{H_0(B-\nu_i)}$ ($n+1\leq i\leq p$), is finite.
\par
Apparently,
$\sigma(A)\sigma_w(B)\cup\sigma_w(A)\sigma(B)\subseteq\sigma(A)\sigma_b(B)\cup\sigma_b(A)\sigma(B)=\sigma_b(\S)$,
see Proposition 4.1(iii).
If $\lambda\notin(\sigma(A)\sigma_w(B)\cup\sigma_w(A)\sigma(B)$),
then, for every $\mu\in\sigma(A)$ and $\nu\in\sigma(B)$ such that
$\mu\nu=\lambda$, $\mu\in\Phi(A)$, $\nu\in\Phi(B)$
and $\ind(A-\mu)=\ind(B-\nu)=0$. Hence $\lambda\in\Phi(\S)$ and
$\ind(\S-\lambda)=0$, i.e., $\lambda\notin\sigma_w(\S)$. Thus
$$\sigma_w(\S)\subseteq\sigma(A)\sigma_w(B)\cup\sigma_w(A)\sigma(B)
\subseteq\sigma(A)\sigma_b(B)\cup\sigma_b(A)\sigma(B)=\sigma_b(\S).$$
The following lemma is the $\sigma_{aw}$ and $\sigma_{ab}$
analogue of this
result.
\begin{lem}\label{lem12}If $A$, $B$ and $\S$ are as above, then
$$\sigma_{aw}(\S)\subseteq\sigma_a(A)\sigma_{aw}(B^*)\cup\sigma_{aw}(A)\sigma_a(B^*)
\subseteq\sigma_a(A)\sigma_{ab}(B^*)\cup\sigma_{ab}(A)\sigma_a(B^*)=\sigma_{ab}(\S).$$
\end{lem}\begin{demo}
The last equality was proved in Proposition 4.3 (iv). The middle
inclusion is a straightforward consequence of the fact that
$\sigma_{aw}(T)\subseteq\sigma_{ab}(T)$ for every operator $T$.
For the first inclusion, let
$\lambda\notin(\sigma_a(A)\sigma_{aw}(B^*)\cup\sigma_{aw}(A)\sigma_a(B^*)$).
Then, for every $\mu_j\in\sigma_a(A)$ and $\nu_j\in\sigma_a(B^*)$
such that $\lambda=\mu_j\nu_j$,  $\mu_j\in\Phi_+(A)$ and
$\nu_j\in\Phi_+(B^*)$. Hence, according to Proposition 4.3(ii),
$\lambda\notin\sigma_{{SF}_+}(\S)$. To complete the proof, it
will be proved that $\ind(\S-\lambda)\leq 0$. Suppose to the
contrary that $\ind(\S-\lambda)>0$ (i.e.,
$\alpha(\S-\lambda)>\beta(\S-\lambda)$). Then
$\lambda\notin\sigma_e(\S)$, and so it follows from the above that
$$\ind(\S-\lambda)=\sum_{j=n+1}^p\ind(A-\mu_j)\dim
H_0(B-\nu_j)-\sum_{j=1}^n\ind (B-\nu_j)\dim H_0(A-\mu_j).$$ Since
$\ind(A-\mu_j)\leq 0$, $\ind(B-\nu_j)\geq 0$ and both
$\dim{H_0(A-\mu_j)}$ and $\dim{H_0(B-\nu_j)}$ are finite,
$\ind(\S-\lambda)\leq 0$, a contradiction. Hence
$\lambda\not\in\sigma_{aw}(\S)$.\end{demo}

\indent The following theorem
gives necessary and sufficient conditions for $\S$ to satisfy gBt
(respectively, $a$-gBt), given that $A$ and $B^*$ satisfy gBt (respectively,
$a$-gBt).

\begin{thm}\label{thm11} Let $\X$ an $\Y$ be two Banach spaces and
consider $A\in B(\X)$ and $B\in B(\Y)$.
\indent {\bf (a).} If $A$ and $B^*$
satisfy gBt, then the following conditions are equivalent:

(i) $\S$ satisfies gBt.

(ii) $\sigma_w(\S)=\sigma(A)\sigma_w(B)\cup\sigma_w(A)\sigma(B)$.

(iii) $A$ has SVEP at points $\mu\in\Phi(A)$ and $B^*$ has SVEP
at points $\nu\in\Phi(B)$ such that
\indent $\mu\nu=\lambda\notin\sigma_w(\S)$.

\

{\bf (b).} If $A$ and $B^*$ satisfy $a$-gBt, then the following
conditions are equivalent:

(I) $\S$ satisfies $a$-gBt.

(II)
$\sigma_{aw}(\S)=\sigma_a(A)\sigma_{aw}(B^*)\cup\sigma_{aw}(A)\sigma_a(B^*)$.

(III) $A$ has SVEP at points $\mu\in\Phi_+(A)$ and $B^*$ has SVEP
at points $\nu\in\Phi_+(B^*)$ such that
$\mu\nu=\lambda\notin\sigma_{aw}(\S)$.\end{thm}

\begin{demo} Since
an operator $T$ satisfies gBt (respectively, $a$-gBt) if and only if $T$
satisfies Bt (respectively, $a$-Bt), \cite{AZ}, it would suffice to prove
that  (i)' $\S$ satisfies Bt $\Longleftrightarrow$ (ii)
$\Longleftrightarrow$ (iii) (resp., (I)' $\S$ satisfies $a$-Bt
$\Longleftrightarrow$ (II) $\Longleftrightarrow$ (III)).

{\bf (a).} The equivalence (i)' $\Longleftrightarrow$ (iii) is
proved in \cite [Theorem 2.1]{D3}; this leaves the equivalence
(i)' $\Longleftrightarrow$ (ii) to be proved. Since $A$ and $B^*$
satisfy Bt, $\sigma_w(A)=\sigma_b(A)$ and
$\sigma_w(B)=\sigma_b(B)$. Hence, if
$\sigma_w(\S)=\sigma(A)\sigma_w(B)\cup\sigma_w(A)\sigma(B)$, then,
according to Proposition 4.2(iii),
$\sigma_b(\S)=\sigma(A)\sigma_b(B^*)\cup\sigma_b(A)\sigma(B^*)
=\sigma(A)\sigma_w(B)\cup\sigma_w(A)\sigma(B)=\sigma_w(\S)$, i.e.,
$\S$ satisfies Bt. (Recall that $\sigma(T)=\sigma(T^*)$ and
$\sigma_b(T)=\sigma_b(T^*)$ for every Banach space operator $T$.)
Conversely, if $\S$ satisfies Bt, then
$\sigma_w(\S)=\sigma_b(\S)=\sigma(A)\sigma_b(B^*)\cup\sigma_b(A)\sigma(B^*)
=\sigma(A)\sigma_w(B)\cup\sigma_w(A)\sigma(B)$.

\

{\bf (b).} (I)' $\Longleftrightarrow$ (II).  Since $A$ and $B^*$
satisfy $a$-Bt, $\sigma_{aw}(A)=\sigma_{ab}(A)$ and
$\sigma_{aw}(B^*)=\sigma_{ab}(B^*)$. Consequently, see
Proposition 4.3(iv) and Lemma \ref{lem12},
$\sigma_{ab}(\S)=\sigma_a(A)\sigma_{ab}(B^*)\cup\sigma_{ab}(A)\sigma_a(B^*)
=\sigma_a(A)\sigma_{aw}(B^*)\cup\sigma_{aw}(A)\sigma_a(B^*)$. Now
if $\S$ satisfies $a$-Bt (so that
$\sigma_{aw}(\S)=\sigma_{ab}(\S)$), then
$\sigma_{aw}(\S)=\sigma_a(A)\sigma_{aw}(B^*)\cup\sigma_{aw}(A)\sigma_a(B^*)$,
and if
$\sigma_{aw}(\S)=\sigma_a(A)\sigma_{aw}(B^*)\cup\sigma_{aw}(A)\sigma_a(B^*)$,
then $\sigma_{aw}(\S)=\sigma_{ab}(\S)$.

(II) $\Longrightarrow$ (III) $\Longrightarrow$ (I)'. For every
factorisation $\lambda=\mu\nu$ of
$\lambda\notin\sigma_{aw}(\S)=\sigma_{ab}(\S)$ such that
$\mu\in\sigma_a(A)$ and  $\nu\in\sigma_a(B^*)$, we have from (II)
that $\mu\in\Phi_+(A)$ with $\asc(A-\mu)<\infty$ and
$\nu\in\Phi_+(B^*)$ with $\asc(B^*-\nu)<\infty$. Since finite
ascent implies SVEP, (II) $\Longrightarrow$ (III) follows.
Suppose now that (III) is satisfied. Let
$\lambda\notin\sigma_{aw}(\S)$. Then
($0\neq$)$\lambda\in\Phi_+(\S)$ and $\ind(\S-\lambda)\leq 0$.
Since, according to Proposition 4.3(ii),
$\sigma_{{SF}_+}(S)=\sigma_a(A)\sigma_{{SF}_+}(B^*)\cup\sigma_{{SF}_+}(A)\sigma_a(B^*)$,
it follows that, for every factorisation $\lambda=\mu\nu$ of
$\lambda$ such that $\mu\in\sigma_a(A)$ and $\nu\in\sigma_a(B^*)$,
$\mu\in\Phi_+(A)$ and $\nu\in\Phi_+(B^*)$. The SVEP hypothesis on
$A$ and $B^*$ now implies that $\asc(A-\mu)$ and $\asc(B^*-\nu)$
are finite, and that $\mu\in iso \hskip.1truecm  \sigma_a(A)$ and $\nu\in \iso\hskip.1truecm  \sigma_a(B^*)$; see Lemma 2.2. Consequently,
$\mu\notin\sigma_{ab}(A)$ and $\nu\notin\sigma_{ab}(B^*)$. This,
by Proposition 4.3(iv), implies that
$\lambda\notin\sigma_{ab}(\S)$, and so
$\sigma_{ab}(\S)\subseteq\sigma_{aw}(\S)$. Since the reverse
inclusion $\sigma_{aw}(T)\subseteq\sigma_{ab}(T)$ holds for every
operator $T$,  see Lemma 4.4, (I)' follows.\end{demo}

\begin{rema}\label{rema11}
{\em The following questions arise naturally from Theorem
\ref{thm11} and its proof. If $A$ and $B^*$ satisfy gBt, then are
the statements (a) $\S$ satisfies gBt  and (b)
$\sigma_{BW}(\S)=\sigma(A)\sigma_{BW}(B)\cup\sigma_{BW}(A)\sigma(B)$
equivalent? (Note that $\sigma_{BW}(T)=\sigma_{BW}(T^*)$ for
every Banach space operator $T$.)  Again, if $A$ and $B^*$ satisfy
$a$-gBt, then are the statements (a)' $\S$ satisfies $a$-gBt and
(b)'
$\sigma_{{SBF}_+^-}(\S)=\sigma_a(A)\sigma_{{SBF}_+^-}(B^*)\cup\sigma_{{SBF}_+^-}(A)\sigma_a(B^*)$
equivalent?}\end{rema} An operator $T$ is said to be polaroid if
points $\lambda\in\iso\hskip.1truecm  \sigma(T)$ are poles of the resolvent, see
\cite{DHJ}. The following lemma proves that the polaroid property
transfers from $A$ and $B$ to $\S$.

\begin{lem}\label{lem13} Let $\X$ and $\Y$ be Banach spaces, and consider $A\in B(\X)$
and $B\in B(\Y)$ two polaroid operators. Then, $\S\in B(\Y,\X)$ is polaroid.\end{lem}

\begin{demo}
Observe that if $\iso\hskip.1truecm  \sigma(A)=\iso \hskip.1truecm  \sigma(B)=\emptyset$, then
$\iso\hskip.1truecm  \sigma(\S)=\emptyset$. If one of the sets $\iso\hskip.1truecm  \sigma(A)$ or
$\iso\hskip.1truecm  \sigma(B)$ is the empty set, say $\iso\hskip.1truecm  \sigma(B)=\emptyset$,
then $\iso \hskip.1truecm  \sigma(\S)\subseteq\{0\}$, $0\in\iso\hskip.1truecm  \sigma(A)$ and
$0\notin\sigma(B)$. Let $\asc(A)=\dsc(A)=d<\infty$. Then
$\asc(L_A)=\dsc(L_A)=d$, see \cite[Theorem 4]{Bo3}. If
$Y={\S}^{d+1}X$ for some $X\in B(\Y,\X)$, then there exists
$Z=AXB\in B(\Y,\X)$ such that $Y={\S}^dZ$, i.e., $\dsc(\S)=d$.
This, since $\S$ has SVEP at $0$, implies that
$\asc(\S)=\dsc(\S)=d$ \cite[Theorem 3.81]{A}. Assume now that
neither of $\iso \hskip.1truecm  \sigma(A)$ and $\iso\hskip.1truecm  \sigma(B)$ is the empty set.
Let $\mu\in\sigma(A)\cap\Pi(A)$ and $\nu\in\sigma(B)\cap\Pi(B)$.
Let $\mu\nu=\lambda$, $\asc(A-\mu)=\dsc(A-\mu)=d_1$,
$\asc(B-\nu)=\dsc(B-\nu)=d_2$ and $d_1+d_2=d$. There exist
decompositions $\X=N((A-\mu)^{d_1})\oplus
(A-\mu)^{d_1}\X={\X}_1\oplus {\X}_2$ and
$\Y=N((B-\nu)^{d_2})\oplus (B-\nu)^{d_2}\Y={\Y}_1\oplus {\Y}_2$
such that $A=A|_{{\X}_1}\oplus A|_{{\X}_2}=A_1\oplus A_2$,
$B=B|_{{\Y}_1}\oplus B|_{{\Y}_2}=B_1\oplus B_2$,
$\sigma(A_1)=\{\mu\}$, $\sigma(B_1)=\{\nu\}$,
$\sigma(A_2)=\sigma(A)\setminus\{\mu\}$,
$\sigma(B_2)=\sigma(B)\setminus\{\nu\}$, $A_1-\mu (=A_1-\mu
I_{{\X}_1})$ is $d_1$--nilpotent and $B_1-\nu$ is
$d_2$--nilpotent. It will be proved that $\dsc(\S-\lambda)\leq
d$, this would then imply that $\lambda\in \Pi(\S)$. Take an $L\in
B(\Y,\X)$ such that $L=(\S-\lambda)^{d+1}M$ for some non--trivial
$M\in B(\Y,\X)$. Let
$M:{\Y}_1\oplus{\Y}_2\rightarrow{\X}_1\oplus{\X}_2$ have the
matrix representation $M=[M_{ij}]_{i,j=1}^2$. Then $$
L=(\S-\lambda)^{d+1}M=\left(\begin{array}{clcr}(L_{A_1}R_{B_1}-\lambda)^{d+1}M_{11}
& (L_{A_1}R_{B_2}-\lambda)^{d+1}M_{12}\\
(L_{A_2}R_{B_1}-\lambda)^{d+1}M_{21} &
(L_{A_2}R_{B_2}-\lambda)^{d+1}M_{22}\end{array}\right).$$ Since
$$L_{A_1}R_{B_1}-\lambda=(L_{A_1}-\mu)R_{B_1}+\mu(R_{B_1}-\nu),$$
where $L_{A_1}-\mu$ is $d_1$--nilpotent and $R_{B_1}-\nu$ is
$d_2$--nilpotent, $L_{A_1}R_{B_1}-\lambda$ is $d$--nilpotent.
Observe that $0\notin\sigma(L_{A_i}R_{B_j}-\lambda)$ for all
$1\leq i, j\leq 2$ such that $i,j\neq 1$. Hence there exist
operators $N_{ij}$ such that
$N_{ij}=(L_{A_i}R_{B_j}-\lambda)M_{ij}$ for all $1\leq i, j\leq
2$; $i,j\neq 1$. Choose $N_{11}\in B({\Y}_1,{\X}_1)$ arbitrarily,
and let $N=[N_{ij}]_{i,j=1}^2$. Then $L=(\S-\lambda)^dN$, i.e.,
$\dsc(\S-\lambda)\leq d$.\end{demo}

An operator $T$ is said to be isoloid (respectively, $a$-isoloid)
if points $\lambda\in\iso\hskip.1truecm  \sigma(T)$ (respectively,
$\lambda\in\iso\hskip.1truecm  \sigma_a(T)$) are eigenvalues of $T$. The isoloid
property transfers from $A$ and $B^*$ to $\S$ \cite{HK}. The
following lemma says that the $a$-isoloid analogue of this result
also holds.
 \begin{lem}\label{lem14} Let $\X$ and $\Y$ be two Banach spaces and consider
$A\in B(\X)$ and $B\in B(\Y)$. If $A$ and $B^*$ are $a$-isoloid,
then $\S$ is $a$-isoloid.\end{lem}\begin{demo} Evidently,
$\iso\hskip.1truecm  \sigma_a(A)=\sigma_a(B^*)=\emptyset\Longrightarrow\iso\hskip.1truecm  \sigma_a(\S)=\emptyset$.
If one of $\iso\hskip.1truecm  \sigma_a(A)$ and $\iso\hskip.1truecm  \sigma_a(B^*)$ is the empty
set, say $\iso\hskip.1truecm  \sigma_a(B^*)=\emptyset$, then
$0\in\iso\hskip.1truecm  \sigma_a(A)$ ($\Longrightarrow 0\in\sigma_p(A)$),
$\iso\hskip.1truecm  \sigma_a(\S)\subseteq\{0\}$ and $0\in\sigma_p(\S)$. If
neither of $\iso\hskip.1truecm  \sigma_a(A)$ and $\iso\hskip.1truecm  \sigma_a(B^*)$ is the empty
set, then every $\lambda\in\iso\hskip.1truecm  \sigma_a(\S)$ has a factoristion
$\lambda=\mu\nu$ such that $\mu\in\iso\hskip.1truecm  \sigma_a(A)$ and
$\nu\in\iso\hskip.1truecm  \sigma_a(B^*)$. $A$ and $B^*$ being $a$-isoloid, it
follows that $\mu\in\sigma_p(A)$ and $\nu\in\sigma_p(B^*)$. Since
$\sigma_p(A)\sigma_p(B^*)\subseteq\sigma_p(\S)$,
$\lambda\in\sigma_p(\S)$.\end{demo}
\indent Recall that an operator $T$
satisfies the generalized Weyl's theorems ($gWt$) if and only if
$\sigma(T)\setminus\sigma_{BW}(T)=\Pi(T)=E(T)$.
\begin{thm}\label{thm12}Let $\X$ and $\Y$ be two Banach spaces and consider
$A\in B(\X)$ and $B\in B(\Y)$.
If $A$ and $B^*$ are isoloid operators which satisfy gWt, then a
necessary and sufficient condition for $\S$ to satisfy gWt is
either (i)
$\sigma_{BW}(\S)=\sigma(A)\sigma_{BW}(B)\cup\sigma_{BW}(A)\sigma(B)$
or (ii)
$\sigma_w(\S)=\sigma(A)\sigma_w(B)\cup\sigma_w(A)\sigma(B)$.\end{thm}\begin{demo}
The hypothesis $A$ and $B^*$ satisfy gWt implies that
$\sigma(A)\setminus\sigma_{BW}(A)=\Pi(A)=E(A)$ and
$\sigma(B)\setminus\sigma_{BW}(B)=\Pi(B)=E(B^*)$.

Assume to start with that (i) is satisfied. Let
$\lambda\in\sigma(\S)\setminus\sigma_{BW}(\S)$. Then, for every
factorisation $\lambda=\mu\nu$ of $\lambda$ such
$\mu\in\sigma(A)$ and $\nu\in\sigma(B)$,
$\mu\notin\sigma_{BW}(A)$ and $\nu\notin\sigma_{BW}(B)$
($=\sigma_{BW}(B^*)$). Consequently, $\mu\in\Pi(A)$ and
$\nu\in\Pi(B)$, and this (by Lemma \ref{lem13}) implies that
$\lambda\in\Pi(\S)$. Thus
$\sigma(\S)\setminus\sigma_{BW}(\S)\subseteq\Pi(\S)$. Since
$\Pi(T)\subseteq\sigma(T)\setminus\sigma_{BW}(T)$ for every
operator $T$, $\sigma(\S)\setminus\sigma_{BW}(\S)=\Pi(\S)$, i.e.
$\S$ satisfies gBt. Observe that $\Pi(\S)\subseteq E(\S)$. If
$\lambda\in E(\S)$, then $\lambda\in\iso\hskip.1truecm  \sigma(\S)$. Since for
every factorisation $\lambda=\mu\nu$ such that $\mu\in\sigma(A)$
and $\nu\in\sigma(B)$, $\mu\in\iso\hskip.1truecm  \sigma(A)$ and
$\nu\in\iso\hskip.1truecm  \sigma(B)$, the isoloid hypothesis on $A$ and $B^*$
implies that $\mu\in E(A)=\Pi(A)$ and $\nu\in E(B^*)=\Pi(B)$.
(Observe that in the case in which $\lambda=0$, and one of
$\iso\hskip.1truecm  \sigma(A)$ and $\iso\hskip.1truecm  \sigma(B)$ is empty, then
$0\in\Pi(\S)$.)  This implies that $\lambda\in\Pi(\S)$. Hence
$E(\S)=\Pi(\S)$, and $\S$ satisfies gWt. Conversely, suppose that
$\S$ satisfies gWt. Then
$\sigma(\S)\setminus\sigma_{BW}(\S)=\Pi(\S)=E(\S)$. If
$\lambda\notin(\sigma(A)\sigma_{BW}(B)\cup\sigma_{BW}(A)\sigma(B)$),
then for every factorisation $\lambda=\mu\nu$ such that
$\mu\in\sigma(A)$ and $\nu\in\sigma(B)$,
$\mu\in\Pi(A)$ and $\nu\in\Pi(B)$. Thus $\lambda\in\Pi(\S)$, which
implies that $\lambda\in\sigma(\S)\setminus\sigma_{BW}(\S)$. Hence
$\sigma_{BW}(\S)\subseteq\sigma(A)\sigma_{BW}(B)\cup\sigma_{BW}(A)\sigma(B)$.
For the reverse inclusion, let
$\lambda\in\sigma(\S)\setminus\sigma_{BW}(\S)$. Then $\lambda\in
E(\S)$, and it follows (see above) that for every factorisation
$\lambda=\mu\nu$ such that $\mu\in\sigma(A)$ and $\nu\in\sigma(B)$,
 $\mu\in\Pi(A)$ and $\nu\in\Pi(B)$. But then
$\mu\notin\sigma_{BW}(A)$ and
$\nu\notin\sigma_{BW}(B)\Longrightarrow \lambda\notin
(\sigma(A)\sigma_{BW}(B)\cup\sigma_{BW}(A)\sigma(B)$).\par To
complete the proof we now consider (ii). Since gWt
$\Longrightarrow$ gBt, see \cite{B6}, the necessity follows from Theorem
\ref{thm11}. Again, Theorem \ref{thm11} implies that if (ii)
holds, then $\S$ satisfies gBt. Now argue as above to prove that
$\S$ satisfies gWt.\end{demo} We say that $\lambda$ is an
$a$-pole of an operator $T$, denoted $\lambda\in\Pi_a(T)$, if
$\lambda\in\iso\hskip.1truecm  \sigma_a(T)$ implies $\asc(T-\lambda)<\infty$ and
the range of $T-\lambda$ is closed.  The operator $T$ is
$a$--polaroid (respectively, left polaroid)  if
$\lambda\in\iso\hskip.1truecm  \sigma_a(T)$ implies $\lambda\in\Pi_a(T)$
(respectively, $\lambda\in\iso\hskip.1truecm  \sigma_a(T)$ implies
$\lambda\in\Pi^a(T)$). Evidently,
$\Pi_a(T)=\sigma_a(T)\setminus\sigma_{ab}(T)$, and $a$-polaroid
implies left polaroid (but not vice--versa). Also, letting
$\Pi_{ao}(T)$ denote finite multiplicity $a$-poles of $T$, $\S$
satisfies $a$-$Bt$ (hence, also $a$-$gBt$) if and only if
$\sigma_a(\S)\setminus\sigma_{aw}(\S)=\Pi_{ao}(\S)$
(equivalently, $\sigma_{aw}(\S)=\sigma_{ab}(\S)$). It is not clear
if the left polaroid (respectively, $a$--polaroid) property
transfers from $A$ and $B^*$ to $\S$, and we have not been
successful in proving an $a$-gWt analogue of Theorem \ref{thm12}.
The following theorem is but a partial analogue of Theorem
\ref{thm12}.

\begin{thm}\label{thm13} Let $\X$ and $\Y$ be two Banach spaces and consider
$A\in B(\X)$ and $B\in B(\Y)$. If $A$ and $B^*$ are $a$-isoloid
operators which satisfy $a$-gWt, then a necessary and sufficient
condition for $\S$ to satisfy $a$-gWt is that
$\sigma_{{SBF}_+^-}(\S)=\sigma_a(A)\sigma_{{SBF}_+^-}(B^*)\cup\sigma_{{SBF}_+^-}(A)\sigma_a(B^*)$
and
$\sigma_{aw}(\S)=\sigma_a(A)\sigma_{aw}(B^*)\cup\sigma_{aw}(A)\sigma_a(B^*)$.\end{thm}\begin{demo}
The hypothesis $A$ and $B^*$ satisfy $a$-gWt  implies that
$\sigma_a(A)\setminus\sigma_{{SBF}_+^-}(A)=\Pi^a(A)=E^a(A)$ and
$\sigma_a(B^*)\setminus\sigma_{{SBF}_+^-}(B^*)=\Pi^a(B^*)=E^a(B^*)$,
and the hypothesis
$\sigma_{aw}(\S)=\sigma_a(A)\sigma_{aw}(B^*)\cup\sigma_{aw}(A)\sigma_a(B^*)$
implies that $\S$ satisfies $a$-gBt (i.e.,
$\sigma_a(\S)\setminus\sigma_{{SBF}_+^-}(\S)=\Pi^a(\S)$); see
Theorem \ref{thm11}({\bf b}). Evidently, $\Pi^a(\S)\subseteq
E^a(\S)$.

{\em Sufficiency.} We already know that $\S$ satisfies $a$-$gBt$.
Let $\lambda\in E^a(\S)$. Then, for every factorisation
$\lambda=\mu\nu$ of $\lambda$ such that $\mu\in\sigma_a(A)$ and
$\nu\in\sigma_a(B^*)$, $\mu\in\iso\hskip.1truecm  \sigma_a(A)$ and
$\nu\in\iso\hskip.1truecm  \sigma_a(B^*)$. The hypothesis $A$ and $B^*$ are
$a$-isoloid thus implies that $\mu\in E^a(A)=\Pi^a(A)$ and
$\nu\in E^a(B^*)=\Pi^a(B^*)$. (If $\lambda=0$, and one of
$\iso\hskip.1truecm  \sigma_a(A)$ and $\iso\hskip.1truecm  \sigma_a(B^*)$ is empty, then
$0\in\Pi^a(\S)$.) Hence
$\mu\in\sigma_a(A)\setminus\sigma_{{SBF}_+^-}(A)$ and
$\nu\in\sigma_a(B^*)\setminus\sigma_{{SBF}_+^-}(B^*)$. This, if
$\sigma_{{SBF}_+^-}(\S)=\sigma_a(A)\sigma_{{SBF}_+^-}(B^*)\cup\sigma_{{SBF}_+^-}(A)\sigma_a(B^*)$,
implies that
$\lambda\in\sigma_a(\S)\setminus\sigma_{{SBF}_+^-}(\S)=\Pi^a(\S)$.
Consequently, $E^a(\S)=\Pi^a(\S)$, and $\S$ satisfies $a$-gWt.\par
{\em Necessity.} If $\S$ satisfies $a$-gWt, then it satisfies
$a$-gBt, and so (see Theorem \ref{thm11}({\bf b}))
$\sigma_{aw}(\S)=\sigma_a(A)\sigma_{aw}(B^*)\cup\sigma_{aw}(A)\sigma_a(B^*)$.
Let $\lambda\in E^a(\S)=
\sigma_a(\S)\setminus\sigma_{{SBF}_+^-}(\S)$. Then, for every
factorisation $\lambda=\mu\nu$ of $\lambda$ such that
$\mu\in\sigma_a(A)$ and $\nu\in\sigma_a(B^*)$, $\mu\in
E^a(A)=\sigma_a(A)\setminus\sigma_{{SBF}_+^-}(A)$ and $\nu\in
E^a(B^*)=\sigma_a(B^*)\setminus\sigma_{{SBF}_+^-}$. (The case in
which one of $\iso\hskip.1truecm  \sigma_a(A)$ and $\iso\hskip.1truecm  \sigma_a(B^*)$ is empty is
dealt with as before.) Consequently,
$\lambda\notin(\sigma_a(A)\sigma_{{SBF}_+^-}(B^*)\cup\sigma_{{SBF}_+^-}(A)\sigma_a(B^*))$,
which implies that
$\sigma_a(A)\sigma_{{SBF}_+^-}(B^*)\cup\sigma_{{SBF}_+^-}(A)\sigma_a(B^*)\subseteq\sigma_{{SBF}_+^-}(\S)$.
For the reverse inclusion, we observe that if
$\lambda\in\sigma_a(\S)$ and
$\lambda\notin(\sigma_a(A)\sigma_{{SBF}_+^-}(B^*)\cup\sigma_{{SBF}_+^-}(A)\sigma_a(B^*)$),
then, for every factorisation $\lambda=\mu\nu$ of $\lambda$ such
that $\mu\in\sigma_a(A)$ and $\nu\in\sigma_a(B^*)$,
$\mu\in\Pi^a(A)=E^a(A)$ and $\nu\in\Pi^a(B^*)=E^a(B^*)$. Thus
$\lambda\in E^a(\S)=\sigma_a(\S)\setminus\sigma_{{SBF}_+^-}(\S)$,
which implies that
$\sigma_{{SBF}_+^-}(\S)\subseteq\sigma_a(A)\sigma_{{SBF}_+^-}(B^*)\cup\sigma_{{SBF}_+^-}(A)\sigma_a(B^*)$.\end{demo}

It is apparent from the proof above that in Theorem \ref{thm13}
one may replace the condition
$\sigma_{aw}(\S)=\sigma_a(A)\sigma_{aw}(B^*)\cup\sigma_{aw}(A)\sigma_a(B^*)$
by the condition that $\S$ satisfies $a$-gBt.\par
\vskip.5truecm
\noindent \bf{Acknowledgements} \rm \par\vskip.3truecm
\indent The authors wish to express their indebtedness to
the referee, for his suggestions
have improved the final version of the present work.\par

\vskip.3truecm
\noindent Enrico Boasso\par
\noindent E-mail address: enrico\_odisseo@yahoo.it\par
\vskip.3truecm
\noindent B. P. Duggal\par
\noindent 8 Redwood Grove, Northfield Avenue,\par
\noindent Ealing, London W5 4SZ,  United Kingdom\par
\noindent E-mail address: duggalbp@gmail.com
 \vskip.3truecm
 \noindent I. H.  Jeon\par
\noindent Seoul National University of Education,\par
\noindent  Seoul 137-742, (South) Korea\par
\noindent E-mail address: jihmath@snue.ac.kr
\end{document}